\newfont{\Bbb}{msbm10 scaled\magstephalf}
\documentclass[12pt]{amsart}
\usepackage{amsmath, amsthm, amscd, amsfonts}
\newcommand{\DD}{{\mathbb D}}
\newcommand{\CC}{{\mathbb C}}

\newcommand{\NN}{{\mathbb N}}
\newcommand{\UU}{{\mathbb U}}

\newcommand{\BB}{{\mathbb B}}
\newcommand{\li}{{\mathcal L}}

\newtheorem{theorem}{Theorem}[section]
\newtheorem{lemma}[theorem]{Lemma}

\newtheorem{corollary}[theorem]{Corollary}

\def\bege{\begin{equation}} \def\ende{\end{equation}}
\def\bsk{\bigskip}   \def\om{\omega}
  \def\ve{\varepsilon} 
\def\pt{\partial}
\def\begr{\begin{eqnarray}} \def\endr{\end{eqnarray}}

\def\bnum{\begin{enumerate}}  \def\enum{\end{enumerate}}

\title[NORM EQUIVALENCE AND COMPOSITION OPERATORS]{NORM EQUIVALENCE AND COMPOSITION
OPERATORS BETWEEN BLOCH/LIPSCHITZ SPACES OF THE UNIT BALL}
\author[Dana Clahane and Stevo Stevi\'c]{Dana D. Clahane \and Stevo Stevi\'c}
\address{
Department of Mathematics\\University of California\\Riverside, CA
92521 USA} \email{dclahane@math.ucr.edu}
\address{Mathematical Institute of the Serbian Academy of Science\\ Knez Mihailova 35/I, 11000 Beograd,
Serbia} \email{sstevic@ptt.yu; sstevo@matf.bg.ac.yu}
\date{}

\begin{document}
\begin{abstract} For $p>0$, let ${\mathcal B}^p(\BB_n)$ and ${\mathcal L}_p(\BB_n)$ respectively
denote the $p$-Bloch and holomorphic $p$-Lipschitz spaces of the
open unit ball $\BB_n$ in $\CC^n$.  It is known that ${\mathcal
B}^p(\BB_n)$ and $\li_{1-p}(\BB_n)$ are equal as sets when
$p\in(0,1)$.  We prove that these spaces are additionally
norm-equivalent, thus extending known results for $n=1$ and the
polydisk. As an application, we generalize work by Madigan on the
disk by investigating boundedness of the composition operator
${\mathfrak C}_\phi$ from ${\mathcal L}_p(\BB_n)$ to ${\mathcal
L}_q(\BB_n)$.
\end{abstract}
\maketitle

\section{Background and Terminology}

Let $n\in\NN$, and suppose that ${\DD}$ is a domain in $\CC^n$.
 Denote the linear space of complex-valued, holomorphic functions on
${\DD}$ by ${\mathcal H}({\DD})$.  If ${\mathcal X}$ is a linear
subspace of ${\mathcal H}({\DD})$ and $\phi:{\DD}\to{\DD}$ is
holomorphic, then one can define the linear operator ${\mathfrak
C}_\phi:{\mathcal X}\rightarrow{\mathcal H}({\DD})$ by ${\mathfrak
C}_\phi(f)=f\circ\phi$ for all $f\in {\mathcal X}$. ${\mathfrak
C}_\phi$ is called the {\em composition operator}$\,$ induced by
$\phi$.

The problem of relating properties of symbols $\phi$ and operators
such as ${\mathfrak C}_\phi$ that are induced by these symbols is of
fundamental importance in concrete operator theory.  However,
efforts to obtain characterizations of self-maps that induce bounded
composition operators on many function spaces have not yielded
completely satisfactory results in the several-variable case,
leaving a wealth of basic, open problems.

In this paper, we try to make progress toward the goal of
characterizing the holomorphic self-maps of the open unit ball
$\BB_n$ in $\CC^n$ that induce bounded composition operators between
holomorphic $p$-Lipschitz spaces ${\mathcal L}_p(\BB_n)$ for $0<p<1$
by translating the problem to $(1-p)$-Bloch spaces ${\mathcal
B}^{1-p}(\BB_n)$ via an auxiliary Hardy/Littlewood-type
norm-equivalence result of potential independent interest.  This
method was also used in \cite{m} for $\BB_1$ and in \cite{csz} for
the unit polydisk $\Delta^n$.

The function-theoretic characterization of analytic self-maps of
$\BB_1$ that induce bounded composition operators on $\li_p(\BB_1)$
for $0<p<1$ is due to K. Madigan \cite{m}, and the case of
$\Delta^n$ was handled in a joint paper by the present authors with
Z. Zhou \cite{csz}, in which a full characterization of the
holomorphic self-maps $\phi$ of $\Delta^n$ that induce bounded
composition operators between ${\mathcal L}_p(\Delta^n)$ and
${\mathcal L}_q(\Delta^n)$, and, more generally, between Bloch
spaces ${\mathcal B}^p(\Delta^n)$ and ${\mathcal B}^q(\Delta^n)$, is
obtained for $p,q\in(0,1)$, along with analogous characterizations
of compact composition operators between these spaces.

Although our main results concerning composition operators, Theorem
\ref{main} and Corollary \ref{lip}, are not full characterizations,
they do generalize Madigan's result for the disk to $\BB_n$; on the
other hand, we obtain a complete Hardy-Littlewood norm-equivalence
result for $p$-Bloch and $(1-p)$-Lipschitz spaces of $\BB_n$ for all
$n\in\NN$. This norm-equivalence result should lead to an eventual
extension to $\BB_n$ of the characterizations of bounded composition
operators established on $\BB_1$ in \cite{m} and on $\Delta^n$ in
\cite{csz}.

Most of our several complex variables notation is adopted from
\cite{r}.  If $z=(z_1,...,z_n)$ and $w=(w_1,...,w_n)$ are points in
$\CC^n$, then we define a complex inner product by $\langle
z,\om\rangle=\sum_{k=1}^n z_k\bar{w_k}$ and put $|z|:=\sqrt{\langle
z, z\rangle}$.  We call $\BB_n:=\{z\in\CC^n:|z|<1\}$ the {\em (open)
unit ball} of $\CC^n$.

Let $p\in (0,\infty).$ The {\em $p$-Bloch space ${\mathcal
B}^p(\BB_n)$} consists of the set of all $f\in {\mathcal H}(\BB_n)$
with the property that there is an $M\geq 0$ such that
$$b(f,z,p):=(1-|z|^2)^p|\nabla f(z)|\leq M\,\,{\text{ for all }}z\in \BB_n$$
${\mathcal B}^p(\BB_n)$ is a Banach space with norm
$||f||_{{\mathcal B}^p}$ given by
\[
||f||_{{\mathcal B}^p}=|f(0)|+\sup_{z\in \BB_n}b(f,z,p).
\]
The {\em little $p$-Bloch space} ${\mathcal B}^p_0(\BB_n)$ is
defined as the closed subspace of ${\mathcal B}^p(\BB_n)$ consisting of the
functions that satisfy
$$\lim_{z\to\pt \BB_n}(1-|z|^2)^p|\nabla f(z)|=0.$$

For $p\in(0,1)$, $\li_p(\BB_n)$ denotes the {\em holomorphic $p$-Lipschitz space}, which is the set of all $f\in
{\mathcal H}(\BB_n)$ such that for some $C>0$,
\begin{equation}\label{heaven}
|f(z)-f(w)|\leq C|z-w|^{p}\,\,\,\,\mbox{for every}\; z,w\in \BB_n.
\end{equation}
These functions extend continuously to ${\overline \BB_n}$ (cf.
\cite[~Lemma~4.4]{csz}).  Therefore, if $A(\BB_n)$ is the ball
algebra \cite[~Ch.~6]{r}, then $${\mathcal L}_p(\BB_n)={\text
{Lip}}_p(\BB_n)\cap A(\BB_n),$$ where Lip$_p(\BB_n)$ is the set of
all $f:\BB_n\rightarrow \CC$ satisfying Equation (\ref{heaven}) for
some $C>0$ and all $z\in{\overline \BB_n}$.  ${\mathcal L}_p(\BB_n)$
is endowed with a complete norm $||\cdot||_{{\mathcal L}_p}$ that is
given by
\begin{equation}\label{haven} ||f||_{{\mathcal L}_p}=|f(0)|+\sup_{z\not= w:
z,w\in{\overline \BB_n}}\left\{\frac{|f(z)-f(w)|}{|z-w|^p}\right\}.
\end{equation}
In Equations (\ref{heaven}) and (\ref{haven}), $\BB_n$ and
${\overline \BB_n}$ are interchangeable, since functions in
${\mathcal L}_p(\BB_n)$ extend continuously to ${\overline
\BB_n}$.  The supremum above is called the {\em Lipschitz
constant} for $f$. As in \cite[~p.~13]{r}, $\sigma$ represents the
unique rotation-invariant positive Borel measure on $\partial
\BB_n$ for which $\sigma(\partial \BB_n)=1$, and for $f\in
L^1(\sigma)$, $C[f]$ denotes Cauchy integral of $f$ on $\BB_n$
(see \cite[~p.~38]{r}).

Let $u\in\partial \BB_n$ and $f\in {\mathcal H}(\BB_n)$. The {\em
directional derivative} of $f$ at $z\in\BB_n$ in the direction of
$u\in\pt\BB_n$ is given by
\[D_uf(z)=\lim_{\lambda\rightarrow 0,\lambda\in\CC}\frac{f(z+\lambda
u)-f(z)}{\lambda}.
\]
Observe that
\begin{equation}\label{gradi}
D_uf(z)=\langle \nabla f(z),{\overline u}\rangle.
\end{equation}We define the partial differential operators
$D_j$ as in \cite[~Ch.~1]{r}.  The radial derivative operator
\cite[~p.~103]{r} in $\CC^n$ will be denoted by ${\mathfrak R}$ and
is linear. Let ${\UU}=\left\{u_1,u_2,\ldots u_n\right\}$ be an
orthonormal basis for the Hilbert space $\CC^n$ with its usual
Euclidean structure.  We define a gradient operator $\nabla^{\UU}$
on ${\mathcal H}({\DD})$ with respect to ${\UU}$ by
\[\nabla^{\UU}f(z)=(D_{u_1}f(z),D_{u_2}f(z),\ldots,
D_{u_n}f(z)),
\]
and we can denote $\nabla^{\UU}$ by $\nabla$ when ${\UU}$ is the
typically ordered standard basis for $\CC^n$.

Let $x$ and $y$ be two positive variable quantities. We write
$x\asymp y$ (and say that $x$ and $y$ are {\em comparable}) if and
only if  $x/y$ is bounded above and below.

\section{Main Results on Composition Operators}
Our norm-equivalence result (Theorem \ref{bigone}) ties our results
concerning ${\mathfrak C}_\phi$ between $p$-Lipschitz spaces of
$\BB_n$ to the following result for general Bloch spaces:
\begin{theorem}\label{main}
Let $p,q\in (0,\infty),$ and suppose that $\phi:\BB_n\rightarrow
\BB_n$ is holomorphic. Then the following statements hold:

(A) If there is an $M\geq 0$ such that for all $z\in\BB_n$ and
$j\in\{1,\ldots,n\}$,
\[
\frac{(1-|z|^2)^q}{(1-|\phi(z)|^2)^p} |\nabla \phi_j(z)|\leq M,
\]

then ${\mathfrak C}_\phi$ is bounded from ${\mathcal B}^p(\BB_n)$
(respectively, ${\mathcal B}^p_0(\BB_n)$) to ${\mathcal
B}^q(\BB_n)$.

(B) If ${\mathfrak C}_\phi$ is bounded from ${\mathcal B}^p(\BB_n)$
(respectively, ${\mathcal B}^p_0(\BB_n)$) to ${\mathcal
B}^q(\BB_n)$, then there is an $M'\geq 0$ such that for all $z\in
\BB_n$ and $u\in\pt\BB_n$,
\[
\frac{(1-|z|^2)^q}{(1-|\langle\phi(z),u\rangle|^2)^p} |\nabla
\langle\phi(z),u\rangle|\leq M'.
\]
\end{theorem}
Theorem \ref{main} above and Corollary \ref{lip} below for $0<p=q<1$
appear in \cite[~Ch.~4]{cl}. It should be pointed out that Theorem
\ref{main}, Part (A) is similar to a statement that is proved in
\cite{zh100}; furthermore, \cite{zh100} contains a result that is in
the same direction as Part (B) of Theorem \ref{main} and that is
proven using different testing functions.  Unlike \cite{zh100},
however, the present paper addresses composition operators between
$\li_p(\BB_n)$ and $\li_q(\BB_n)$ and the coincidence and
norm-equivalence of ${\mathcal B}^{1-p}(\BB_n)$ and ${\mathcal
L}_p(\BB_n)$, respectively.

It is natural to consider the application of corresponding
``little-oh" arguments to obtain a compactness result analogous to
Theorem \ref{main}, in which ``bounded" is replaced by ``compact"
and the limit of the left hand side of each inequality in the
statement is taken as $|\phi(z)|\rightarrow 1^-$, with inequality
replaced by equality to $0$.  However, in the case that $p\in(0,1)$,
${\mathcal B}^p(\BB_n)$ is the same as and norm-equivalent to
${\mathcal L}_{1-p}(\BB_n)$, whose compact composition operators are
known (by a result due to J. H. Shapiro) to be generated precisely
by holomorphic self-maps $\phi$ of $\BB_n$ with supremum norm
strictly less than $1$ (see \cite[~Ch.~4]{cm}).

The following corollary follows from Theorems \ref{main} and
\ref{bigone} and extends the main result of \cite{m}:
\begin{corollary}\label{lip}
Let $\,\,p,q\in(0,1)$, and suppose that $\phi:\BB_n\rightarrow
\BB_n$ is holomorphic.  Then the following statements hold:

(A) If there is an $M\geq 0$ such that
\[
\frac{(1-|z|^2)^{1-q}}{(1-|\phi(z)|^2)^{1-p}} | \nabla
\phi_j(z)|\leq M,
\]
for all $j\in\{1,2\ldots,n\}$ and $z\in\BB_n$, then ${\mathfrak
C}_\phi$ is a bounded operator from ${\mathcal L}_p(\BB_n)$ to
${\mathcal L}_q(\BB_n)$.

(B) If ${\mathfrak C}_\phi$ is a bounded operator from ${\mathcal
L}_p(\BB_n)$ to ${\mathcal L}_q(\BB_n)$, then there is an $M'\geq 0$
such that for all $z\in \BB_n$ and $u\in\pt\BB_n$,
\[
\frac{(1-|z|^2)^{1-q}}{(1-|\langle\phi(z),u\rangle|^2)^{1-p}}
|\nabla \langle\phi(z),u\rangle|\leq M'.
\]
\end{corollary}
Choosing $n=1$, $p=q\in(0,1)$, and $u=1$ in Corollary \ref{lip}
leads to the following result, which is due to K. Madigan \cite{m}:
\begin{theorem}\label{Madig}
Let $\,\,0<p<1$, and suppose that $\phi$ is an analytic self-map of
$\BB_1$.  Then ${\mathfrak C}_\phi$ is bounded on ${\mathcal
L}_p(\BB_1)$ if and only if
\[
\sup_{z\in \BB_1}\left\{\left(\frac{1-|z|^2}{1-|\phi(z)|}\right)^
{1-p}|\phi'(z)|\right\}<\infty.
\]
\end{theorem}
\section{Norm Equivalence of ${\mathcal L}_p(\BB_n)$ and ${\mathcal B}^{1-p}(\BB_n)$ for $0<p<1$.}

To generalize Theorem \ref{Madig} to $\BB_n$, we need Theorem
\ref{bigone}, which is the ball analogue of the following result for
the disk (Lemma 2 in \cite{m}). The first statement in Theorem
\ref{H-L} can be derived from a classical theorem of
Hardy/Littlewood for $n=1$ (see \cite{hl}, \cite[~p.~74]{d}, and
\cite[~p.~176]{cm}).
\begin{theorem}\label{H-L}
Let $0<p<1$.  If $f:\BB_1\rightarrow \CC$ is analytic, then $f\in $
${\mathcal L}_p(\BB_1)$ if and only if
\[
|f'(z)|=O\left(\frac{1}{1-|z|^2}\right)^{1-p}.
\]
Furthermore, the Lipschitz constant of $f$ and the quantity
\[
\sup_{z\in \BB_1}\{(1-|z|^2)^{1-p}|f'(z)|\}
\]
are comparable as $f$ varies through ${\mathcal L}_p(\BB_1)$.
\end{theorem}
We remark that the polydisk version of Theorem \ref{H-L} is stated
and proved in \cite{csz}.  However, the argument used there cannot
be applied to $\BB_n$, so we need a different approach for that
domain.  We will proceed by listing some lemmas, which together
eseentially form the norm equivalence Theorem \ref{bigone}.

For $0<p<1$, we can define a norm $||f||_{{\mathcal
B}^{1-p}}^{\mathfrak R}$ on $\li_p(\BB_n)$ by
\[
   ||f||_{{\mathcal
B}^{1-p}}^{\mathfrak R}=|f(0)|+\sup_{z\in \BB_n}
\{(1-|z|^2)^{1-p}|({\mathfrak R}f)(z)|\}.\]  The following lemma is
part of our norm equivalence result, Theorem \ref{bigone}:
\begin{lemma}\label{dont2}
Suppose that $0<p<1$.   Furthermore, there is a $C_p\geq 0$ such
that for all $f\in \li_p(\BB_n)$,
\[
||f||_{{\mathcal B}^{1-p}}^{\mathfrak R}\leq C_p||f||_{{\mathcal
L}_p}.
\]
\end{lemma}
\begin{proof}
The proof of the first statement is standard and left to the reader.
 Since functions in ${\mathcal L}_p(\BB_n)$ extend continuously to
${\overline \BB_n}$, they are automatically in $L^1(\sigma)$
\cite[~Remark,~p.~107]{r} and since the quotients of these functions
and their $\li_p$-norms satisfy \cite[~Equation~(1),~p.~107]{r}, the
second statement is obtained from \cite[~Theorem~6.4.9]{r}.
\end{proof}
The following lemma is also a portion of Theorem \ref{bigone}:
\begin{lemma}\label{grr}
If $p\in(0,1)$, then ${\mathcal B}^{1-p}(\BB_n)\subset\li_p(\BB_n)$,
and
\[||f||_{{\mathcal L}_p}\leq (2+2p^{-1})||f||_{{\mathcal
B}^{1-p}}\,\,\mbox{ for all }\,f\in {\mathcal B}^{1-p}(\BB_n).
\]
\end{lemma}
\begin{proof} Suppose that $f\in{\mathcal B}^{1-p}(\BB_n)$.  If $f=0$ then $f\in{\mathcal L}_p(\BB_n)$
trivially, so assume henceforward that $f\not=0$.
 A well-known result \cite[~Ch.~6]{r} applied to $f/||f||_{{\mathcal B}^{1-p}}$ implies that for all $z,w\in B_n$,
\[
\frac{1}{||f||_{{\mathcal B}^{1-p}}}|f(z)-f(w)|\leq
(1+2p^{-1})|z-w|^p,
\]
from which the first statement of the lemma follows.  Moreover,
\begin{eqnarray*}
||f||_{{\mathcal L}_p}&=&
|f(0)|+\sup_{z,w\in\BB_n:z\not=w}\frac{|f(z)-f(w)|}{|z-w|^{p}}\\
&\leq &|f(0)|+(1+2p^{-1})||f||_{{\mathcal
B}^{1-p}}\\
&\leq &(2+2p^{-1})||f||_{{\mathcal B}^{1-p}}.
\end{eqnarray*}
\end{proof}

The following fact also constitutes part of Theorem \ref{bigone}:
\begin{lemma}\label{st}
If $p>0$, then $f\in{\mathcal B}^p(\BB_n)$ if and only if there
exists $M\geq 0$ such that $|({\mathfrak R}f)(z)|(1-|z|^2)^p\leq M$
for all $z\in\BB_n$.  Also, $||\cdot||_{{\mathcal B}^p}^{\mathfrak
R}$ given by $||f||_{{\mathcal B}^p}^{\mathfrak
R}:=|f(0)|+\sup_{z\in\BB_n}|({\mathfrak R}f)(z)|(1-|z|^2)^p$ is a
norm on ${\mathcal B}^p(\BB_n).$  If $p\in (0,1]$, then there is a
$C_p\geq 0$ such that $||f||_{{\mathcal B}^p}\leq
C_p||f||_{{\mathcal B}^p}^{\mathfrak R}$ for all $f\in{\mathcal
B}^p(\BB_n).$
\end{lemma}
\begin{proof}
For a proof of the first statement, see \cite[Proposition 1]{yo}.
The second statement follows from subsequent applications of the
first statement in Lemma \ref{grr} and Lemma \ref{dont2}. To prove
the final statement, we use the weighted Bergman projection $P_s$
with kernel $K_s$ and the map $L_s$ defined on $P_s[L^\infty(B_n)]$
by
\[
(L_sg)(z)=(s+1)^{-1}(1-|z|^2)\left[(n+s+1)g(z)+({\mathfrak
R}g)(z)\right]\mbox{ for all }z\in\BB_n,
\]
where $s\in\CC$ satisfies Re\,\,$s>-1$ (see \cite{c}). By
\cite[Corollary 13]{c}, we have that $P_s\circ L_s$ is the identity
on ${\mathcal B}^1(\BB_n)$ for all such values of $s$. In
particular, $P_0\circ L_0$ is the identity on ${\mathcal
B}^p(\BB_n)$, since this set is contained in ${\mathcal
B}^1(\BB_n)$.  Note that the assumption $p\in(0,1]$ is used here.

We then obtain that there is a $C\geq 0$ such that for all $z\in
B_n$ and $f\in{\mathcal B}^1(\BB_n)$,
\begin{eqnarray*}
f(z)&=& (P_0\circ L_0)(f)(z)\\
&=&C \int_{\BB_n} (1-|w|^2) K_0(z,w)\biggl[(n+1) f(w)+{\mathfrak R}
f(w)\biggr]  dV(w). \end{eqnarray*} Hence, there is a $C'\geq 0$
such that for all $f\in{\mathcal B}^p(\BB_n)$ and $z\in \BB_n$,
\begr |\nabla f(z)| &\le& C' \int_{\BB_n} (1-|w|^2) |\nabla K_0(z,w)|\,  |f(w)|  dV(w)\nonumber\\
&&+C' \int_{\BB_n} (1-|w|^2) |\nabla K_0(z,w)|\, |{\mathfrak R}
f(w)| dV(w).\nonumber\endr

\noindent Let $\ve\in (1-p,1)$.  Subsequent applications of the
above inequality, \cite[~Lemma~2]{s}, and \cite[Theorem 1.4.10]{r}
imply that there are non-negative constants $C''$ and $C'''$ such
that for all $z\in B_n$ and $f\in {\mathcal B}^p(\BB_n)$, the
following chain of inequalities holds: \begr|\nabla f(z)|
&\le& C' \int_{\BB_n} \frac{(1-|w|^2) |w|}{|1-\langle z, w\rangle|^{n+2}} |f(w)| dV(w)\nonumber\\
&&+C' \int_{\BB_n} \frac{(1-|w|^2) |w|}{|1-\langle z, w\rangle|^{n+2}}|{\mathfrak R} f(w)|  dV(w)\nonumber\\
&\le&C'' ||f||_{{\mathcal
B}^p}^{\mathfrak R}\int_{\BB_n} \frac{(1-|w|^2)^{\ve}}{|1-\langle z, w \rangle|^{n+2}}  dV(w)\nonumber\\
&&+C'' ||f||_{{\mathcal
B}^p}^{\mathfrak R}\int_{\BB_n} \frac{(1-|w|^2)^{1-p}}{|1-\langle z, w \rangle|^{n+2}}  dV(w)\le\nonumber\\
&\le& C'''||f||_{{\mathcal B}^p}^{\mathfrak
R}\frac1{(1-|z|)^{1-\ve}}+ C'''||f||_{{\mathcal B}^p}^{\mathfrak R}\frac1{(1-|z|)^{p}}.\nonumber\endr\\
It follows that for all $f\in{\mathcal B}^p(\BB_n)$ and $z\in\BB_n$,
$$(1-|z|^2)^p|\nabla f(z)|\leq 2^{p+1}C'''||f||_{{\mathcal B}^p}^{\mathfrak R}.$$
The final statement in the lemma now follows from the above
statement and an application of \cite[~Lemma~2]{s} at $z=0$.
\end{proof}
Next, we state and prove this section's main result, the analogue of
Theorem \ref{H-L} for $\BB_n$.  We emphasize that while the
statement of equality in the theorem is known and can be obtained,
for example, from \cite{zhu}, the norm equivalence portion requires
additional work that includes the previous lemmas and the proof
below. Furthermore, neither this result nor its proof has appeared
previously in any literature that is known to the authors, though it
seems to be part of the folklore.  The proof of this rather
fundamental theorem seems to be non-trivial and worthy of recording.
\begin{theorem}\label{bigone}
If $0<p<1$, then ${\mathcal B}^{1-p}(\BB_n)=\li_p(\BB_n)$;
furthermore,
\[||f||_{{\mathcal B}^{1-p}}\asymp ||f||_{{\mathcal
B}^{1-p}}^{\mathfrak R}\asymp ||f||_{{\mathcal L}_p}
\]
as $f$ varies through ${\mathcal L}_p(\BB_n)$.
\end{theorem}

\begin{proof} The first statement is known, since ${\mathcal L}_p(\BB_n)=A(\BB_n)\cap $Lip$_\alpha(\BB_n)$
(see \cite[~Ch.~6]{r}), which is set-theoretically equal to
${\mathcal B}^{1-p}(\BB_n)$ (see \cite{yo}). By Lemma \ref{st}, it
follows that there is a $C_p\geq 0$ such that for all
$f\in\li_p(\BB_n)$, $||f||_{{\mathcal B}^{1-p}}\leq
C_p||f||_{{\mathcal B}^{1-p}}^{\mathfrak R}$.  It follows from Lemma
\ref{dont2} that there is a $C_p'\geq 0$ such that for all
$f\in\li_p(\BB_n)$, $||f||_{{\mathcal B}^{1-p}}\leq
C_p||f||_{{\mathcal B}^{1-p}}^{\mathfrak R}\leq
C_pC_p'||f||_{{\mathcal L}_p}$, which is less than or equal to
$C_pC_p'(2+2p^{-1})||f||_{{\mathcal B}^{1-p}}$ by Lemma \ref{grr}.
The second statement in Theorem \ref{bigone} follows.
\end{proof}

\section{Proof of Theorem \ref{main}}\label{simp}

In the proof of Theorem \ref{main}, part (B), we will use part of
the following lemma, which is obtained by straightforward estimates
involving Equation (\ref{gradi}) (see \cite[~Ch.~4]{cl}):
\begin{lemma}\label{ugrad} Let $f\in {\mathcal H}({\DD}),$ where ${\DD}$ is an open subset of $\CC^n$, and suppose that ${\UU}$ is an
orthonormal basis for $\CC^n.$ Then for all $z\in{\DD}$,
\[
|\nabla^{\UU}f(z)|\asymp|\nabla f(z)|.
\]
\end{lemma}
We are now ready to prove Theorem \ref{main}.\bsk

{\it Proof of Theorem \ref{main}.  (A)} Suppose that for some $M\geq
0$,
\begin{eqnarray}
\quad\,\,\,\frac{(1-|z|^2)^q}{(1-|\phi(z)|^2)^p} |\nabla
\phi_j(z)|\leq M\label{condi}\mbox{ for all } \,z\in
\BB_n,\,j\in\{1,2,\ldots ,n\}.
\end{eqnarray}
If $z\in\BB_n$ and $F(z)=(1-|z|^2)^q|\nabla({\mathfrak C}_\phi
f)(z)|$.  Then we have that
\begin{eqnarray}
F(z)&=&(1-|z|^2)^q\sqrt{\sum_{i=1}^n\left|D_i(f\circ\phi)(z)\right|^2}\nonumber\\
&\leq&
(1-|z|^2)^q\sum_{i=1}^n \left|D_i(f\circ\phi)(z)\right|\nonumber\\
&\leq &(1-|z|^2)^qn\sum_{j=1}^n \left|\nabla f(\phi(z))\right|\left|\nabla \phi_j(z)\right|\nonumber\\
&=&n \left|\nabla f(\phi(z))\right|(1-|\phi(z)|^2)^p \frac{(1-|z|^2)^q}{(1-|\phi(z)|^2)^p}\sum_{j=1}^n \left|\nabla \phi_j(z)\right|\nonumber\\
&\leq&n \sup_{w\in \BB_n}\left\{\left|\nabla
f(w)\right|(1-|w|^2)^p\right\}
\sum_{j=1}^n \frac{(1-|z|^2)^q}{(1-|\phi(z)|^2)^p}\left|\nabla \phi_j(z)\right|\nonumber\\
&\leq& n||f||_{{\mathcal B}^p}nM\label{hy},
\end{eqnarray}
by Inequality (\ref{condi}).  It follows that $||{\mathfrak C}_\phi
f||_{{\mathcal B}^q}\leq (1+n^2M)||f||_{{\mathcal B}^p}$ for every
$f\in {\mathcal B}^p(\BB_n)$, thus completing the proof of Theorem
\ref{main}, Part (A).

{\it (B).} We proceed by modifying the argument given in
\cite[~p.~187-188]{cm} for $n=1$.  For $a\in \BB_n$, define
$f_a:\BB_n\rightarrow \CC$ to be function that vanishes at $0$ and
is the antiderivative of $\psi_a:\BB_n\rightarrow\CC$ given by
$\psi_a(t)=(1-\bar a t)^{-p}$.  Let $w\in \BB_n$ and $u\in\partial
\BB_n$.  Define $F_{w,u}:\BB_n\rightarrow\CC$ by
\[
F_{w,u}(z)=f_{\langle w,u\rangle}(\langle z,u\rangle).
\]
Define $\phi_u:\BB_n\rightarrow \BB_1$ by $\phi_u(z)=\langle
\phi(z),u\rangle$.  Let $u^{(1)}:=u$, and choose
$u^{(2)},u^{(3)},\ldots ,u^{(n)}$ so that
${\UU}=\{u^{(1)},u^{(2)},u^{(3)},\ldots ,u^{(n)}\}$ is an
orthonormal basis for $\CC^n$.  For all $z\in \BB_n$ and
\,$j\in\{2,3,\ldots,n\}$, we have that
\begin{eqnarray}
D_{u^{(j)}}F_{w,u}(z)&=&\lim_{\lambda\rightarrow 0}\frac{F_{w,u}(z+\lambda u^{(j)})-F_{w,u}(z)}{\lambda}\nonumber\\
&=&\lim_{\lambda\rightarrow 0}\frac{f_{\langle w,u^{(1)}\rangle}(\langle z+\lambda u^{(j)},u^{(1)}\rangle)- f_{\langle w,u^{(1)}
\rangle}(\langle z,u^{(1)}\rangle)}{\lambda}\nonumber\\
&=&0\label{tes}.
\end{eqnarray}
On the other hand, for every $z\in \BB_n$,
\begin{eqnarray}
D_{u^{(1)}}F_{w,u}(z)
&=&\lim_{\lambda\rightarrow 0}\frac{F_{w,u}(z+\lambda u)-F_{w,u}(z)}{\lambda}\nonumber\\
&=&\lim_{\lambda\rightarrow 0}\frac{f_{\langle w,u\rangle}(\langle
z,u\rangle+\lambda) -
f_{\langle w,u\rangle}(\langle z,u\rangle)}{\lambda}\nonumber\\
&=& \psi_{\langle w,u\rangle}(\langle z,u\rangle)\label{tes1}.
\end{eqnarray}
From Equations (\ref{tes}) and (\ref{tes1}), it follows that
\begin{equation}\label{tes3}
|\nabla^{\UU}F_{w,u}(z)|=|\psi_{\langle w,u\rangle}(\langle
z,u\rangle)|= \left|1-\overline{\langle w,u\rangle}\langle
z,u\rangle\right|^{-p}.
\end{equation}
We observe that the quantity above is bounded when $u$ is fixed.
This fact and Lemma \ref{ugrad} together imply that $F_{w,u}\in
{\mathcal B}_0^p(\BB_n)$. Also, we have
\[F_{w,u}(0)=f_{\langle w,u\rangle}(\langle 0,u\rangle)=f_{\langle w,u\rangle}(0)=0.\]
Furthermore, by Lemma \ref{ugrad}, we have that
\begin{eqnarray}
\sup_{z\in \BB_n}(1-|z|^2)^{p}|\nabla F_{w,u}(z)|
&=&\sup_{z\in \BB_n}(1-|z|^2)^{p}|\nabla^{\UU} F_{w,u}(z)|\nonumber\\
&=&\sup_{z\in \BB_n}(1-|z|^2)^{p}\left|1-{\overline{\langle
w,u\rangle}}\langle z,u\rangle\right|^{-p}\label{van1}.
\end{eqnarray}
Note that
\[
\left|1-{\overline{\langle w,u\rangle}}\langle
z,u\rangle|\right|^{-p}\leq (1-|z|)^{-p}\leq
\frac{2^{p}}{(1-|z|^2)^{p}}.
\]
It follows that Quantity (\ref{van1}) is less than or equal to
$2^{p}$.  Hence, $F_{w,u}\in {\mathcal B}^p(\BB_n)$ for every $w\in
\BB_n$ and $u\in\pt\BB_n$; moreover, the set
\[\{||F_{w,u}||_{{\mathcal
B}^p}: u\in\pt\BB_n,\,\,w\in \BB_n\}
\]
is bounded.  This fact and the hypothesis together imply that there
exist $C$ and $M\geq 0$ such that for every $w\in \BB_n$ and
$u\in\pt\BB_n$,
\[
||F_{w,u}\circ\phi||_{{\mathcal B}^q}\leq C||F_{w,u}||_{{\mathcal
B}^p}\leq CM.
\]
Therefore, we obtain that
\begin{equation}\label{tric}
\sup_{u\in\pt\BB_n,\,\,z,w\in \BB_n}\left\{|\nabla (f_{\langle
w,u\rangle}\circ\phi_{u})(z)|(1-|z|^2)^{q}\right\}\leq CM.
\end{equation}
Now for each $j\in\{1,2,\ldots,n\}$, we have that
\begin{eqnarray*}
D_{j} (f_{\langle w,u\rangle}\circ\phi_{u})(z)
&=& f'_{\langle w,u\rangle}(\langle \phi(z),u\rangle)D_j\langle \phi(z),u\rangle\\
&=& \left(1-{\overline{\langle w,u\rangle}}\langle
\phi(z),u\rangle\right)^{-p} D_j\langle \phi(z),u\rangle.
\end{eqnarray*}
It follows that
\begin{equation}\label{glass}
\nabla(f_{\langle w,u\rangle}\circ\phi_{u})(z)
=\left(1-{\overline{\langle w,u\rangle}}\langle
\phi(z),u\rangle\right)^{-p} \nabla\langle \phi(z),u\rangle
\end{equation}
Using Equation (\ref{glass}), we can rewrite Equation (\ref{tric}) as
\[
\sup_{u\in\pt\BB_n,\,\,z,w\in \BB_n}
\frac{(1-|z|^2)^{q}}{\left|1-{\overline{\langle w,u\rangle}}\langle
\phi(z),u\rangle\right|^{p}}|\nabla\langle \phi(z),u\rangle|\leq CM.
\]
In particular, we have that
\begin{equation}\label{endy}
\sup_{u\in \pt \BB_n,\,\, z\in
\BB_n}\frac{(1-|z|^2)^{q}}{(1-|\langle
\phi(z),u\rangle|)^{p}}|\nabla\langle \phi(z),u\rangle|\leq CM,
\end{equation}
from which the statement of Theorem \ref{main}, Part (B) follows.\bsk

By restricting the values of $u$, one obtains various necessary
conditions for compactness of $C_\phi$ from Part (B) of Theorem
\ref{main}.  Two such conditions are listed in Corollary
\ref{062005} below.  We point out that the boundedness of Quantity
(\ref{pet}) below when ${\mathfrak C}_\phi$ is bounded from
${\mathcal B}^p(\BB_n)$ to ${\mathcal B}^q(\BB_n)$ is a result given
by Zhou in \cite{zh100}.
\begin{corollary}\label{062005} Let $p,q>0$.  If ${\mathfrak C}_\phi$ is a bounded operator from ${\mathcal B}^p(\BB_n)$ (respectively,
${\mathcal B}^p_0(\BB_n))$ to ${\mathcal B}^q(\BB_n)$, then there is
an $M\geq 0$ such that the following statements hold:

(i) For all $z\in \BB_n$ with $\phi(z)\not=0$, we have that
\begin{eqnarray}
\frac{(1-|z|^2)^q}{(1-|\phi(z)|^2)^p}\frac{|J_\phi(z)^T\cdot\phi(z)|}{|\phi(z)|}\leq
M.\label{pet}
\end{eqnarray}

(ii) For all $z\in\BB_n$ and $j\in\{1,2,\ldots,n\}$,
\begin{equation}\label{ebob}
\frac{(1-|z|^2)^q}{(1-|\phi_j(z)|^2)^p}|\nabla \phi_j(z)|\leq M.
\end{equation}
\end{corollary}

\begin{proof} Putting $u:={\overline{\phi(z)}}/|\phi(z)|$ in Theorem \ref{main}, Part (B),
one obtains that Quantity (\ref{pet}) is no larger than some $M'\geq
0$ for all $z\in B_n$ such that $\phi(z)\not=0$. Successively
replacing $u\in\pt\BB_n$ in Theorem \ref{main}, Part (B) by the
typically ordered standard basis elements $e_j$ of $\CC^n$ for
$j=1,2,\ldots,n$, we see that the left side of Inequality
(\ref{ebob}) is no larger than some $M''\geq 0$, so that we can
choose $M:=\max(M',M'')$.
\end{proof}

\section{Acknowledgements}\nonumber
The authors would like to thank W. Wogen for kindly pointing out
that an earlier, coordinate-dependent version of Theorem \ref{main},
Part (B) could be improved to its current coordinate-free form.


\begin{thebibliography}{99}
\bibitem[Cho]{c} B. R. Choe. {\em Projections, the weighted Bergman spaces, and the Bloch space},
{\em Proc. Amer. Math. Soc.} {\bf 108} (1990), 127-136.
\bibitem[Cl]{cl} D. D. Clahane. {\em Composition operators on
holomorphic function spaces of several compex variables}, Ph.D.
Thesis, University of California, Irvine, 2000.
\bibitem[CSZ]{csz} D. D. Clahane/S. Stevi\'c/Z. Zhou. {\em Composition operators on general Bloch spaces of the
polydisk}, preprint.
\bibitem[CoMac]{cm} C. C. Cowen and B. D. MacCluer.  {\em Composition Operators
on Spaces of Analytic Functions}, CRC Press, Boca Roton, 1995.
\bibitem[D]{d} P. L. Duren. {\em Theory of $H^p$ Spaces}, Academic Press, New York, 1970.
\bibitem[HL]{hl} G. H. Hardy/J. E. Littlewood. {\em Some properties of fractional integrals II},
{\it Math. Z.} {\bf 34} (1932), 403-439.
\bibitem[Mad]{m} K. Madigan. {\em Composition operators
on analytic Lipschitz spaces}, {\it Proc. Amer. Math. Soc.} {\bf
119} (2) (1993), 471-480.
\bibitem[Ru]{r} W. Rudin.  {\em Function Theory in the Unit Ball of ${\CC}^n$,} Springer-Verlag, New York, 1980.
\bibitem[S]{s} S. Stevi\'c. {\em On an integral operator on the unit
ball in ${\bf C}^n$}, J. Inequal. Appl. {\bf 1} (2005), 81-88.
\bibitem[YO]{yo} W. Yang/C. Ouyang. {\em Exact location of $\alpha$-Bloch spaces
in $L^p_a$ and $H^p$ of a complex unit ball}, Rocky Mountain J.
Math. {\bf 30} (2000), no. 3, 1151-1169.
\bibitem[Zho]{zh100} Z. Zhou. {\em Composition operators between
$p$-Bloch space and $q$-Bloch space in the unit ball}, Progress in
Nat. Sci. {\bf 13} (2003), no. 3, 233-236.
\bibitem[Zhu]{zhu} K. Zhu. {\em Spaces of Holomorphic Functions in
the Unit Ball}, Springer, New York, 2005.
\end{thebibliography}
\end{document}